\documentclass{amsart}
\usepackage{natbib}

\theoremstyle{definition}
\newtheorem{definition}{Definition}[section]

\theoremstyle{plain}
\newtheorem{theorem}[definition]{Theorem}
\newtheorem{lemma}[definition]{Lemma}
\newtheorem{corollary}[definition]{Corollary}
\newtheorem{conjecture}[definition]{Conjecture}
\newtheorem{proposition}[definition]{Proposition}
\newtheorem{remark}[definition]{Remark}

\def\zvm#1#2#3#4{\left(\begin{array}{cc}#1&#2\\#3&#4\end{array}\right)}

\title[Moufang loops of order $64$]
{Toward the Classification of Moufang Loops of Order $64$}

\author{Petr Vojt\v{e}chovsk\'y}

\address{Department of Mathematics, University of Denver, 2360 S Gaylord St,
Denver, 80208, Colorado, U.S.A.}

\email{petr@math.du.edu}

\begin{document}

\begin{abstract} We show how to obtain all nonassociative Moufang loops of
order less than $64$ and $4262$ nonassociative Moufang loops of order $64$ in a
unified way. We conjecture that there are no other nonassociative Moufang loops
of order $64$. The main idea of the computer search is to modify precisely one
quarter of the multiplication table in a certain way, previously applied to
small $2$-groups.
\end{abstract}



\maketitle

\section{Introduction}

\noindent A set $Q$ with one binary operation is a \emph{quasigroup} if the
equation $xy = z$ has a unique solution in $Q$ whenever two of the three
elements $x$, $y$, $z\in Q$ are specified. \emph{Loop} is a quasigroup with a
neutral element $1$ satisfying $1x=x1=x$ for every $x$. \emph{Moufang loops}
are loops in which any of the (equivalent) \emph{Moufang identities}
\begin{align}
    ((xy)x)z &= x(y(xz)),\label{Eq:M1}\tag{M1}\\
    x(y(zy)) &= ((xy)z)y,\label{Eq:M2}\tag{M2}\\
    (xy)(zx) &= x((yz)x),\label{Eq:M3}\tag{M3}\\
    (xy)(zx) &= (x(yz))x\label{Eq:M4}\tag{M4}
\end{align}
holds. It was shown recently \cite{PhVoMonthly} that, in an analogy to
groups, any set with one binary operation, neutral element and two-sided
inverses satisfying either \eqref{Eq:M1} or \eqref{Eq:M2} is already a
Moufang loop.

Moufang loops are certainly the most studied loops. They arise naturally in
algebra (as the multiplicative loop of octonions \cite{SpVe2000},
\cite{CoSm2003}), and in projective geometry (Moufang planes \cite{TiWe2003}),
for example.

Although Moufang loops are generally nonassociative, they retain many
properties of groups that---borrowing a phrase from \cite[p.\ 7]{Co2001}---we
know and love. For instance: (i) every $x$ is accompanied by its two-sided
inverse $x^{-1}$ such that $xx^{-1}=x^{-1}x=1$, (ii) any two elements generate
a subgroup (this property is called \emph{diassociativity}), (iii) in finite
Moufang loops, the order of an element divides the order of the loop, and, as
is believed to be shown recently in \cite{GrZa2003}, the order of a subloop
divides the order of the loop.

On the other hand, many essential tools of group theory are not available for
Moufang loops. The lack of associativity makes presentations very awkward and
hard to calculate, and permutation representations in the usual sense
impossible.

It is therefore no surprise that the classification of Moufang loops of order
$n$ is completed only up to and including $n=63$ \cite{Ch1978},
\cite{GoMaRa1999}. Several ingenious constructions, described in detail in
\cite{GoMaRa1999}, are needed to obtain all the loops.

In this paper, we introduce a class of Moufang loops that includes all
nonassociative Moufang loops of order less than $64$, and $4262$ nonassociative
Moufang loops of order $64$ (compare this with the $267$ groups of order $64$).
We conjecture that there are no other nonassociative Moufang loops of order
$64$.

The class is obtained by a computer program based on an idea of Dr\'apal. It
takes only a few minutes to obtain the Moufang loops of order less than $64$,
and about $2$ weeks to obtain $4262$ Moufang loops of order $64$ (using a PC
with $2$ GHz processor).

Thanks to this algorithm, small Moufang loops can now be stored in a uniform
and very efficient way (about $4$ bytes of data are needed for a Moufang loop
of order $64$). They are available via the GAP \cite{GAP} package LOOPS
\cite{LOOPS} written by G.~Nagy and the present author. Great care was taken to
comply with the naming conventions introduced in \cite{GoMaRa1999}.

Unfortunately, there is no guarantee that the algorithm found all
nonassociative Moufang loops of order $64$, and, in fact, it is not clear how
this question could be answered easily. Nevertheless, it appears to be a
definite step toward the classification of small Moufang loops, especially
small Moufang $2$-loops.

\subsection{Organization of this paper.}
It is known that a finite Moufang loop has order $p^n$ if and only if it has
exponent $p^m$, for some prime $p$ and integers $n$, $m$. This fact is recalled
and newly proved in Section \ref{Sc:2ML}.

Dr\'apal's cyclic and dihedral constructions are described in Section
\ref{Sc:CD}, where we also summarize some results of these constructions
obtained in an earlier paper \cite{DrVo2004}.

The computer search always starts with a single Moufang loop, referred to as a
\emph{seed}. We use the so-called loops $M(G,2)$ (due to Chein) as seeds. The
definition and properties of the loops $M(G,2)$ can be found in Section
\ref{Sc:MG2}.

The computer search is outlined in Section \ref{Sc:Outline}, where we also
present the results in a tabular form. The reader who is only interested in the
outcome of the search will understand it fully at that point and does not have
to read further.

The algorithm is discussed in detail in Section \ref{Sc:Algorithm}.

Several nontrivial theoretical results were needed to make the algorithm
sufficiently fast. These are collected and proved in Section \ref{Sc:SpeedUp}.
We pay attention especially to the isomorphism problem for (Moufang) loops.

Section \ref{Sc:GAP} contains detailed instructions on how to obtain and use
the GAP package LOOPS.

The paper closes with a section devoted to conjectures and open problems.

\section{Moufang $2$-loops}\label{Sc:2ML}

\noindent A loop is said to be \emph{power associative} if the power $x^n$ is
well-defined for every element $x$ and a positive integer $n$. Moufang loops
are power associative, by diassociativity.

Let $p$ be a prime. We say that a power associative loop has \emph{exponent}
$p^r$ if the order of every element of $L$ divides $p^r$. Finite power
associative loops of exponent $p^r$, for some $r$, are called \emph{$p$-loops}.

One of the fundamental facts of group theory is that a finite group has
exponent $p^r$ if and only if it is of order $p^s$. This certainly does not
generalize to $p$-loops. It is easy to construct by hand a loop of order $5$
and exponent $2$, for instance. Another well-known example is the smallest
nonassociative Steiner loop of order $10$ and exponent $2$ \cite{CoRo1999}.

This has the unfortunate consequence that the two natural definitions of a
$p$-loop are not equivalent, yet they appear side by side in the literature.
Since we deal predominantly with Moufang loops of order $64=2^6$ here, let us
first make sure that all is well for Moufang loops. The following proposition
was first proved by Glauberman \cite{Gl} for odd $p$, and by Glauberman and
Wright \cite{GlWr} for $p=2$. We offer a short proof that relies on the
classification of finite simple Moufang loops, and hence on the
classification of finite simple groups. The original proofs of Glauberman and
Wright do not require the classification.

Recall that a subloop $H$ of a loop $L$ is \emph{normal} in $L$ if $aH=Ha$,
$a(bH)=(ab)H$, and $(aH)b=a(Hb)$ holds for every $a$, $b\in L$. Given
elements $x$, $y$, $z$ of a loop $L$, the \emph{associator} $[x,y,z]\in L$ is
defined by $(xy)z=(x(yz))[x,y,z]$. The \emph{associator subloop} $A(L)$ of
$L$ is the subloop of $L$ generated by all associators $[x,y,z]$. The
\emph{nucleus} $N(L)$ of $L$ consists of all elements $x\in L$ such that
$[x,y,z]=[y,x,z]=[y,z,x]=1$ for every $y$, $z\in L$. Finally, the
\emph{center} $Z(L)$ is the subloop $\{x\in N(L);\; xy=yx$ for every $y\in
L\}$.

\begin{proposition} Let $M$ be a finite Moufang loop and $p$ a prime. Then $M$
is of exponent $p^r$ for some $r$ if and only if $M$ has order $p^s$ for some
$s$.
\end{proposition}
\begin{proof}
Let $|M|=p^s$ and let $x\in M$. As is well-known (cf.\ \cite[p.\ 13]{Pf1990}),
the order of $x$ divides the order of $M$. In particular, $M$ is of exponent
$p^s$.

Conversely, suppose that $M$ is of exponent $p^s$. Assume, for a contradiction,
that $|M|$ is not a power of $p$, and that $s$ is as small as possible.

If $M$ is not simple, it possesses a nontrivial normal subloop $L$. Then
$|M|=|L|\cdot |M/L|$. Both $L$, $M/L$ are Moufang loops of exponent a power of
$p$. Since $|L|<|M|$ and $|M/L|<|M|$, the orders of $L$ and $M/L$ are powers of
$p$, by the induction hypothesis. Then $|M|$ is a power of $p$, too.

We complete the proof by showing that there is no nonassociative finite simple
Moufang loop of exponent $p^s$.

Liebeck classified all nonassociative finite simple Moufang loops in
\cite{Li1987}. It turns out that there is exactly one nonassociative finite
simple Moufang loop $M^*(q)$ for every finite field $GF(q)$. The loops $M^*(q)$
are obtained as follows (see \cite{Pa1956}, \cite{VoJGT} for more details):

Let $F=GF(q)$. Consider the \emph{Zorn vector matrices}
\begin{equation}\label{Eq:GenericZVM}
    \zvm{a}{\alpha}{\beta}{b},
\end{equation}
where $a$, $b\in F$, and $\alpha$, $\beta\in F^3$. The matrices are
multiplied according to the Zorn multiplication formula
\begin{displaymath}
    \zvm{a}{\alpha}{\beta}{b} \zvm{c}{\gamma}{\delta}{d}
    = \zvm{ac+\alpha\cdot\delta}{a\gamma+\alpha d-\beta\times\delta}{\beta c
    +b\delta + \alpha\times\gamma}{\beta\cdot\gamma +bd},
\end{displaymath}
where $\alpha\cdot\beta$ (resp.\ $\alpha\times\beta$) is the dot product (resp.
cross product) of $\alpha$ and $\beta$.

Let $M(q)$ consist of all matrices \eqref{Eq:GenericZVM} with
$ab-\alpha\cdot\beta=1$. Then $M^*(q)=M(q)/Z(M(q))$. Note that the group
$PSL(2,q)$ embeds into $M^*(q)$ via
\begin{displaymath}
    \zvm{a}{b}{c}{d}\mapsto \zvm{a}{(b,0,0)}{(c,0,0)}{d},
\end{displaymath}
since all cross products vanish when two such vector matrices are multiplied.
Since no $PSL(2,q)$ is a $p$-group, we are done.
\end{proof}


\section{The Cyclic and Dihedral Constructions}\label{Sc:CD}

\noindent While working on the problem of Hamming distances of groups
\cite{Dr1992}, Dr\'apal discovered two constructions that modify exactly one
quarter of the multiplication table of a group and yield another group, often
with a different center and thus not isomorphic to the original group. It is
known \cite{Dr2000} that two $2$-groups whose multiplication tables (with
rows and columns labelled in the same way) coincide in more than three
quarters of the cells must be isomorphic. Hence, the two constructions
exemplify a minimal change in a $2$-group (in the sense of multiplication
tables) that yields a nonisomorphic group.

Let us first give a brief description of the constructions and then talk about
their power. Note that the constructions work for Moufang loops, too. The
generalization from groups to Moufang loops was carried through in
\cite{DrVo2004}.

\subsection{Modular arithmetic}
For a positive integer $m$, let $M=\{-m+1,\dots,m\}$. Define $\sigma:\mathbb
Z\to \{-1,0,1\}$ by
\begin{displaymath}
    \sigma(i)=\left\{\begin{array}{ll}
        0,&i\in M,\\
        1,&i>m,\\
        -1,&i<-m+1.
    \end{array}\right.
\end{displaymath}
Then $\sigma$ can be used to describe addition $\oplus$ and subtraction
$\ominus$ modulo $M$. Namely, for $i$, $j\in M$, we have $i\oplus j =
i+j-2m\sigma(i+j)$, $i\ominus j = i-j-2m\sigma(i-j)$.

\subsection{The cyclic construction}

Let $G$ be a Moufang loop with normal subloop $S$ such that $G/S$ is a cyclic
group of order $2m$. Let $\alpha$ be a generator of $G/S$. Then for every
$x\in G$ there is a unique $i\in M=\{-m+1,\dots,m\}$ such that
$x\in\alpha^i\subseteq G$. Let $h\in Z(G)\cap S$. Define a new multiplication
$*$ on $G$ by
\begin{displaymath}
    x*y = xyh^{\sigma(i+j)},
\end{displaymath}
where $x\in\alpha^i$, $y\in\alpha^j$, $i\in M$, $j\in M$. Note that no
parentheses are needed in the expression $xyh^{\sigma(i+j)}$ because $h\in
N(G)$.

As is shown in \cite{DrVo2004}, the resulting loop $(G,*)$ is a Moufang loop.
Adopting the notation of \cite{Dr2003}, the loop $(G,*)$ will also be denoted
by $G[S,\alpha,h]$ or $G[\alpha,h]$.

\subsection{The dihedral construction}

Let $G$ be a Moufang loop with normal subloop $S$ such that $G/S$ is a
dihedral group of order $4m$, $m\ge 1$. Let $\beta$, $\gamma$ be involutions
of $G/S$ such that $\alpha=\beta\gamma$ is of order $2m$. Pick $e\in\beta$,
$f\in\gamma$. Then for every element $x\in G$ there are uniquely determined
integers $i$, $j\in M$ such that $x\in\alpha^i\cup e\alpha^i$ and
$x\in\alpha^j\cup \alpha^jf$. Let $G_0=\bigcup_{i\in M}\alpha_i \le G$. Let
$h\in N(G)\cap Z(G_0)\cap S$. Define a new multiplication $*$ on $G$ by
\begin{displaymath}
    x*y = xyh^{(-1)^r\sigma(i+j)},
\end{displaymath}
where $x\in\alpha^i\cup e\alpha^i$, $y\in\alpha^j\cup \alpha^jf$, $i\in M$,
$j\in M$, and $r\in\{0,1\}$ is equal to $0$ if and only if $y\in G_0$.

As is shown in \cite{DrVo2004}, the resulting loop $(G,*)$ is a Moufang loop,
and will also be denoted by $G[S,\beta,\gamma,h]$ or $G[\beta,\gamma,h]$.

\subsection{Pictorial representation of the constructions}

The reader might get a better feel for the constructions when considering the
effect of the constructions on the multiplication table of $G$. The diagrams in
Figure \ref{Fg:Diags} indicate the changes to the multiplication table caused
by the cyclic construction for $m=4$ (left) and by the dihedral construction
for $m=2$ (right). Each square represents $(n/|S|)^2$ elements of $G$. The
multiplication table of $(G,*)$ differs from the multiplication table of
$(G,\cdot)$ according to the symbol in the square: no symbol $\Rightarrow$ no
change, ``$+$'' $\Rightarrow$ multiply every entry by $h$, ``$-$''
$\Rightarrow$ multiply every entry by $h^{-1}$. Viewed in this way, the
constructions get a more combinatorial flavor.

\begin{figure}
\begin{footnotesize}
\begin{displaymath}
\setlength\arraycolsep{2pt}
\begin{array}{c||c|c|c|c|c|c|c|c|}
    *&\alpha^0&\alpha^1&\alpha^2&\alpha^3&\alpha^4&\alpha^{-3}&\alpha^{-2}&\alpha^{-1}\\
    \hline
    \hline
    \alpha^0&&&&&&&&\\ \hline
    \alpha^1&&&&&+&&&\\ \hline
    \alpha^2&&&&+&+&&&\\ \hline
    \alpha^3&&&+&+&+&&&\\ \hline
    \alpha^4&&+&+&+&+&&&\\ \hline
    \alpha^{-3}&&&&&&-&-&-\\ \hline
    \alpha^{-2}&&&&&&-&-&\\ \hline
    \alpha^{-1}&&&&&&-&&\\ \hline
\end{array}\quad\quad
\begin{array}{c||c|c|c|c|c|c|c|c|}
    *&\alpha^0&\alpha^0f&\alpha^1&\alpha^1f&\alpha^2&\alpha^2f&
    \alpha^{-1}&\alpha^{-1}f\\
    \hline
    \hline
    \alpha^0&&&&&&&&\\ \hline
    e\alpha^0&&&&&&&&\\ \hline
    \alpha^1&&&&&+&-&&\\ \hline
    e\alpha^1&&&&&+&-&&\\ \hline
    \alpha^2&&&+&-&+&-&&\\ \hline
    e\alpha^2&&&+&-&+&-&&\\ \hline
    \alpha^{-1}&&&&&&&-&+\\ \hline
    e\alpha^{-1}&&&&&&&-&+\\ \hline
\end{array}
\end{displaymath}
\end{footnotesize}
\caption{Pictorial representation of the constructions.} \label{Fg:Diags}
\end{figure}

\subsection{Invariants of the constructions}

The essential properties of the constructions are summarized in the following
theorem:

\begin{theorem}[Theorem 6.3, Theorem 6.4 \cite{DrVo2004}]\label{Th:OldMain}
Let $(G,\cdot)$ be a Moufang loop and let $(G,*)$ be obtained from $(G,\cdot)$
by the cyclic or the dihedral construction. Then:
\begin{enumerate}
\item[(i)] $(G,*)$ is a Moufang loop,

\item[(ii)] $(G,*)$ is a group if and only if $(G,\cdot)$ is,

\item[(iii)] the associators of $(G,\cdot)$, $(G,*)$ are in $S$, and thus the
associator subloops of $(G,\cdot)$ and $(G,*)$ coincide as loops,

\item[(iv)] the nuclei of $(G,\cdot)$ and $(G,*)$ coincide as sets,

\item[(v)] the constructions are reversible, i.e., $(G,\cdot)$ is obtained from
$(G,*)$ by the cyclic or the dihedral construction with some parameters.
\end{enumerate}
\end{theorem}

\emph{Extra loops} are loops satisfying the identity $x(y(zx)) = ((xy)z)x$.
Extra loops are precisely Moufang loops with all squares in the nucleus
\cite[Corollary 2]{ChRo1972}. The constructions preserve this property of
Moufang loops:

\begin{lemma}\label{Lm:Extra} Let $(G,\cdot)$ be an
extra loop and let $(G,*)$ be obtained from $(G,\cdot)$ by the cyclic or the
dihedral construction. Then $(G,*)$ is an extra loop.
\end{lemma}
\begin{proof}
It suffices to show that $x*x\in N(G,*)$ for every $x\in G$. This follows
immediately from Theorem \ref{Th:OldMain}(iv), since $x*x=x^2h^\varepsilon$
for some $\varepsilon$, and $x^2\in N(G,\cdot)$, $h\in N(G,\cdot)$.
\end{proof}

\subsection{Using the constructions}
A good question is whether the constructions are powerful enough to produce
many $2$-groups from a single group. Given two groups $G$, $H$, let us call $H$
a \emph{modification} of $G$ if there is an integer $n$ and groups $G=K_0$,
$K_1$, $\dots$, $K_{n-1}$, $K_n\cong H$ such that $K_{i+1}$ is obtained from
$K_i$ by one of the two constructions, for $0\le i<n$. For a group $G$, let
$\mathcal M(G)$ denote the set of all modifications of $G$. We will call $G$
the \emph{seed} of $\mathcal M(G)$.

Since the constructions are reversible, every element of $\mathcal M(G)$ is in
fact a seed of $\mathcal M(G)$. The most optimistic plan is therefore to show
that given any group $G$ of order $2^m$, $\mathcal M(G)$ comprises all groups
of order $2^m$. This indeed happens for $n=2^m\le 32$. (This was noticed by the
present author for $n=8$ in \cite{VoMSThesis}, and by Dr\'apal and Zhukavets
for $n=16$, $32$ in \cite{DrZh2003}.)

There are, of course, other, much faster means of generating $2$-groups (cf.
the manual of GAP \cite{GAP} or the survey paper \cite{BEO}), however, none of
the group-theoretical methods applies to Moufang loops.

Theorem \ref{Th:OldMain} claims that the nuclei and associator subloops are
invariant under the constructions. A quick glance into the classification of
small Moufang loops \cite{GoMaRa1999} reveals that some Moufang loops of
order 32 have nucleus of size $2$, others of size $4$. Hence no single
nonassociative Moufang loop of order 32 can possible yield all other Moufang
loops of that order by a repeated application of the two constructions,
shattering our most optimistic plan outlined above. We need more seeds.


\section{Seeds for the computer search}\label{Sc:MG2}

\noindent There is a class of nonassociative Moufang loops, first defined by
Chein \cite{Ch1978}, that is well understood. Let $G$ be a group of order $n$,
and let $u$ be a new element. Define multiplication $\circ$ on $G\cup Gu$ by
\begin{displaymath}
    g\circ h = gh,\quad g\circ hu = (hg)u,\quad gu\circ h = (gh^{-1})u,
    \quad gu\circ hu = hg^{-1},
\end{displaymath}
where $g$, $h\in G$. The resulting loop $(G\cup Gu,\circ)=M(G,2)$ is a Moufang
loop. It is nonassociative if and only if $G$ is nonabelian.

We are going to show that $M(G,2)$ is isomorphic to $M(H,2)$ if and only if
$G$ is isomorphic to $H$. Thus, we will obtain as many nonassociative Moufang
loops of order $2n$ as there are nonabelian groups of order $n$. Proposition
\ref{Pr:MG2} is probably well known, but since we were unable to find a
reference, we give a proof here.

For a finite power-associative loop $L$ and a positive integer $i$, let
\begin{displaymath}
    s_i(L)=|\{x\in L;\; |x|=i\}|.
\end{displaymath}
We call $s(L)=(s_1(L)$, $s_2(L)$, $\dots)$ the \emph{order statistic} of $L$.
The following Lemma shows why black-box recognition of finite abelian groups is
not hard in principle.

\begin{lemma}\label{Lm:AOS}
A finite abelian group is determined uniquely by its order statistic.
\end{lemma}
\begin{proof}
Let $A$ be a finite abelian group. For a prime $p$, let $A_{(p)}=\{x\in
A;x^{(p^i)}=1$ for some $i\}$ be the $p$-primary component of $A$. Then
$s_{p^k}(A_{(p)}) = s_{p^k}(A)$. Thus, $s(A)$ determines $s(A_{(p)})$, and it
suffices to prove the lemma for all finite abelian $p$-groups $A$.

Let $m$ be the largest integer with $s_{p^m}(A)>0$. Then $A=B\times C$, where
$C$ is a cyclic group of order $p^m$, and $B$ is a finite abelian $p$-group.
As $A=B\times C$ is a direct product, we have
\begin{displaymath}
    s_{p^u}(A)=s_{p^u}(B)\cdot\Bigl(|C|-\sum_{v>u}s_{p^v}(C)\Bigr) +
        \Bigl(|B|-\sum_{v>u}s_{p^v}(B)\Bigr)\cdot s_{p^u}(C) -
        s_{p^u}(B)s_{p^u}(C).
\end{displaymath}
Since the order statistics of $A$ and $C$ are known, the order statistic of $B$
can be calculated, starting with $s_{p^m}(B)$. We are done by induction on
$|A|$.
\end{proof}

\begin{proposition}\label{Pr:MG2}
Assume that $G$, $H$ are two finite groups. Then $G\cong H$ if and only if
$M(G,2)\cong M(H,2)$.
\end{proposition}
\begin{proof}
Only one implication is nontrivial. Assume that $M(G,2)\cong M(H,2)$. Then we
can consider $H$ to be a subgroup of $M(G,2)$. By \cite[Lemma 3.11]{ChGo2003}
or by \cite[Subsection 4.2]{VoThesis}, either $H=G$ (and we are done), or there
is a subgroup $A$ of $G$ such that $H=M(A,2)$. Since $H$ is associative, $A$ is
abelian. Similarly, either $G=H$ (and we are done), or there is an abelian
group $B$ such that $G\cong M(B,2)$.

The order statistic of a group $K$ and the order statistic of the associated
loop $M(K,2)$ can be reconstructed from each other, because the coset $Ku$
consists of involutions. Being isomorphic, the loops $M(G,2)=M(M(B,2),2)$ and
$M(H,2)=M(M(A,2),2)$ have identical order statistics. Thus the abelian groups
$B$, $A$ have identical order statistics, and are isomorphic by Lemma
\ref{Lm:AOS}. Then $G=M(B,2)$, $H=M(A,2)$ are isomorphic, too.
\end{proof}


\section{Notation and results of the computer search}\label{Sc:Outline}

\noindent From now on, whenever we say \emph{Moufang loop} we mean a
\emph{nonassociative Moufang loop}.

Given a seed (Moufang loop) $M$, we can calculate the class of Moufang loops
$\mathcal M(M)$, collecting only one loop of each isomorphism type.

Thanks to Section \ref{Sc:MG2}, we have plenty of seeds with which to start
the computer search. It turns out that all Moufang $2$-loops of order less
than $64$ are obtained from the seeds $M(G,2)$, and only four more seeds (see
below) are needed in addition to the loops $M(G,2)$ to obtain all Moufang
loops of order less than $64$.

The results of the search can be found in Table \ref{Tb:Main}. Here is how to
read Table \ref{Tb:Main}.

Under \texttt{class}, we give the name of the class $\mathcal M(M)$ of Moufang
loops. The names are systematic if the seed is of the form $M(G,2)$, and
\emph{ad hoc} in the $4$ remaining cases.

When the seed of order $2n$ is of the form $M(G,2)$, then $G$ is a nonabelian
group of order $n$. (Table \ref{Tb:NonAb} gives the number of nonabelian groups
of order $1\le n\le 32$ with orders for which no nonabelian group exists
omitted.) Each such group is identified uniquely in GAP (version 4.3). If it is
cataloged as the $m$th nonabelian group of order $n$ in GAP, it can be obtained
by the GAP command \texttt{AllGroups(n, IsCommutative, false)[m]}, for
instance. Accordingly, we use the name $2n:m$ for the corresponding class of
Moufang loops. (Warning: Since we cannot guarantee that the GAP libraries of
groups will not change in the future, the reader should note the version of GAP
carefully.)

When the seed of order $2n$ is not of the form $M(G,2)$, we denote the class by
$2n:xm$, as in $36:x1$.

Under \texttt{|nucleus|}, we give the size of the nucleus of all loops in the
class.

Under \texttt{assoc.\ subloop}, we give the isomorphism type of the associator
subloop of all loops in the class, using standard group-theoretical notation.
Hence, $C_m$ denotes the cyclic group of order $m$, $Q_8$ denotes the
quaternion group of order $8$, and $A_4$ denotes the alternating group of order
$12$.

Under \texttt{seed(s)}, we list the seed that was used to generate the class.
When an integer $m$ is listed, the seed is the loop $M(G,2)$ where $G$ is the
$m$th nonabelian group of order $n$. When several integers are listed, then all
corresponding loops $M(G,2)$ can be used as seeds, but only the first one was
actually used in the search. In the remaining cases $2n:xm$, we give the seed
explicitly by referring to smaller Moufang loops. Here,
\texttt{MoufangLoop(n,m)} denotes the $m$th Moufang loop of order $n$, as
cataloged in \cite{GoMaRa1999} and in the package LOOPS.

Under \texttt{extra?} we specify if all loops in the class are extra loops
(yes), or if all loops in the class are not extra loops (no). No other
scenarios can occur by Lemma \ref{Lm:Extra}.

Under \texttt{|class|} we specify the number of nonisomorphic loops forming the
class.

\begin{table}
\caption{Classes of nonassociative Moufang loops obtained by the cyclic and
dihedral constructions from the indicated seeds. All nonassociative Moufang
loops of order less than $64$ are accounted for in this table.}\label{Tb:Main}
\begin{displaymath}
    \begin{array}{c|c|c|c|c|r}
    \text{class}&\text{$|$nucleus$|$}&\text{assoc.\ subloop}&
    \text{seed(s)}&
    \text{extra?}&\text{$|$class$|$}\\
    \hline \hline
    12:01&1&C_3&1&\text{no}&1\\
    \hline
    16:01&2&C_2&1,2&\text{yes}&5\\
    \hline
    20:01&1&C_5&1&\text{no}&1\\
    \hline
    24:01&2&C_3&1,3&\text{no}&4\\
    24:02&1&C_2\times C_2&2&\text{no}&1\\
    \hline
    28:01&1&C_7&1&\text{no}&1\\
    \hline
    32:01&4&C_2&\text{1--3, 7--9}&\text{yes}&60\\
    32:04&2&C_4&\text{4-6}&\text{no}&11\\
    \hline
    36:01&1&C_9&1&\text{no}&1\\
    36:02&3&C_3&2&\text{no}&1\\
    36:03&1&C_3\times C_3&3&\text{no}&1\\
    36:x1&3&C_3&\mathrm{MoufangLoop}(12,1)\times C_3&\text{no}&1\\
    \hline
    40:01&2&C_5&\text{1, 3}&\text{no}&4\\
    40:02&1&C_5&2&\text{no}&1\\
    \hline
    42:01&1&C_7&1&\text{no}&1\\
    \hline
    44:01&1&C_{11}&1&\text{no}&1\\
    \hline
    48:01&4&C_3&\text{1, 4, 5, 12}&\text{no}&19\\
    48:02&2&Q_8&2&\text{no}&2\\
    48:03&2&C_6&\text{3, 5, 7}&\text{no}&11\\
    48:08&6&C_2&\text{8, 9}&\text{yes}&11\\
    48:10&1&A_4&10&\text{no}&1\\
    48:11&2&C_2\times C_2&11&\text{no}&2\\
    48:x1&6&C_2&\mathrm{MoufangLoop}(16,4)\times C_3&\text{yes}&5\\
    \hline
    52:01&1&C_{13}&1&\text{no}&1\\
    \hline
    54:01&3&C_3&1&\text{no}&1\\
    54:02&3&C_3&2&\text{no}&1\\
    \hline
    56:01&2&C_7&\text{1, 2}&\text{no}&4\\
    \hline
    60:01&5&C_3&1&\text{no}&1\\
    60:02&3&C_5&2&\text{no}&1\\
    60:03&1&C_{15}&3&\text{no}&1\\
    60:x1&3&C_5&\mathrm{MoufangLoop}(20,1)\times C_3&\text{no}&1\\
    60:x2&5&C_3&\mathrm{MoufangLoop}(12,1)\times C_5&\text{no}&1\\
    \hline
    64:01&8&C_2&\text{1--3, 10, 14, 18--22, 32, 33}&\text{yes}&1316\\
    64:04&2&C_2\times C_2&\text{4--6}&\text{no}&18\\
    64:07&4&C_4&\text{7--9, 11--13, 34--37}&\text{no}&214\\
    64:15&2&C_8&\text{15--17}&\text{no}&11\\
    64:23&4&C_2\times C_2&\text{23--31}&\text{yes}&2612\\
    64:38&2&C_4&\text{38, 39}&\text{no}&44\\
    64:43&2&C_2&\text{43, 44}&\text{yes}&47\\
    \hline
    \hline
    \end{array}
\end{displaymath}
\end{table}

\begin{table}
\caption{Number of isomorphism classes of nonabelian groups of order $1\le n\le
32$.}\label{Tb:NonAb}
\begin{small}
\begin{displaymath}
\begin{array}{r|cccccccccccccccc}
    \text{order}&6&8&10&12&14&16&18&20&21&22&24&26&27&28&30&32\\
    \text{nonab.\ groups}&1&2&1&3&1&9&3&3&1&1&12&1&2&2&3&44
\end{array}
\end{displaymath}
\end{small}
\end{table}

\subsection{What the results indicate.}
As we have already mentioned, both constructions preserve the nucleus (as a
set) and the associator subloop (as a loop). Let us therefore say that the
\emph{parameter} of a seed $M$ is the triple $(|M|$, $|N(M)|$, isomorphism
class of $A(M))$.

With one exception (classes $54:01$, $54:02$), two seeds $M(G,2)$ are in the
same class of loops if and only if their parameters coincide.

More importantly, the seeds $M(G,2)$ generate all Moufang $2$-loops of order
less than $64$, and all but $4$ classes of Moufang loops of order less than
$64$. The four exceptional cases are all generated by seeds of the form
$M\times C_{2k+1}$, where $M$ is a Moufang loop of smaller order.

Table \ref{Tb:Main} accounts for all Moufang loops of order less than $63$,
according to the classification \cite{GoMaRa1999}.

\begin{remark}\label{Rm:1} It is known that the $267$ groups of order $64$ split into
two classes $($of size $261$ and $6)$ with respect to the modifications. We
have checked that none of the $6$ groups in the second class is of the form
$M(G,2)$, where $G$ is a group of order $32$.
\end{remark}


\section{The Algorithm}\label{Sc:Algorithm}

\noindent This section describes the main steps of the algorithm used to
calculate the class $\mathcal M(M)$ from a seed $M$.

\subsection{Platform} All calculations were implemented in GAP version 4.3
for Windows, using the package LOOPS. The search ran for about $2$ weeks on a
PC with a $2$GHz processor.

\subsection{Input.} A Moufang loop $M$ (seed), flagged as unexplored.

\subsection{Output.} The class of Moufang loops $\mathcal M(M)$ (with one
Moufang loop for every isomorphism type) together with data that describes how
to build all loops of $\mathcal M(M)$ from the seed $M$.

\subsection{Main cycle.} Let $L$ be the first unexplored loop. If there is
none, the search is over. Otherwise:
\begin{enumerate}
\item[(i)] determine all normal subloops $S$ of $L$ such that $L/S$ is cyclic
of even order or dihedral of doubly even order,

\item[(ii)] for every normal subloop $S$ of $L$, determine all admissible
parameters of the constructions of Section \ref{Sc:CD} (e.g., in the cyclic
case, find all pairs $(\alpha,h)$ where $\alpha$ is a generator of $L/S$ and
$h\in S\cap Z(L)$),

\item[(iii)] using the parameters found in step (ii), construct the
modifications $(L,*)$ from $L$,

\item[(iv)] store those newly found loops $(L,*)$ that are not isomorphic to
any of the previously found loops; flag them as unexplored,

\item[(v)] flag $L$ as explored.
\end{enumerate}


\section{Speeding up the algorithm}\label{Sc:SpeedUp}

\noindent The steps (i), (ii) and (iv) are expensive, especially step (iv).
We describe in this section how to speed up (ii) and (iv). Many additional
improvements of programming character were incorporated into the algorithm
but we do not mention them here.

\subsection{Speeding up step (ii)}
The problem with step (ii) is that there are typically very many parameters
$S$, $\alpha$, $\beta$, $\gamma$, $h$ that can be used to modify the loop $L$
into $(L,*)$. Since we are only interested in the isomorphism type of the
resulting loop $(L,*)$, we would like to know which parameters yield isomorphic
loops. This topic has been studied for groups in \cite{Dr2003}. For example, it
is proved in \cite{Dr2003} that the cyclic modification $G[S,\alpha,h]$ is
independent of the generator $\alpha$ of $S$, in the sense that for two
generators $\alpha$, $\alpha'$ of $S$ and $h\in S\cap Z(L)$ there is $h'\in
S\cap Z(L)$ such that $G[S,\alpha,h]$ is isomorphic to $G[S,\alpha',h']$. Such
an observation speeds up the search substantially, since a cyclic group of
order $n$ contains $\varphi(n)$ generators, where $\varphi$ (the Euler
function) counts the number of positive integers relatively prime to $n$.

Unfortunately, it is by no means easy to generalize the results of
\cite{Dr2003} into the nonassociative case. (In fact, it is often impossible,
for we have found counterexamples to some generalizations of \cite{Dr2003}.)
What follows is a generalization of the above result (independence of
$\alpha$ in the cyclic construction) for a class of Moufang loops with the
associator subloop contained in the center. By \cite{KiKu2004}, all extra
$2$-loops $L$ of order less than $512$ satisfy $A(L)\subseteq Z(L)$. Table
\ref{Tb:Main} shows that the two largest classes of Moufang loops of order
$64$ consist of extra loops.

We follow the reasoning of \cite{Dr2003}, often word for word. The proofs had
to be expanded substantially when diassociativity did not apply.

\begin{lemma}\label{Lm:Aux} Let $\sigma$ be as in Section $\ref{Sc:CD}$.
For every $i$, $j\in M=\{-m+1,\dots,m\}$, we have:
\begin{enumerate}
\item[(i)]  $\sigma(i+j)+\sigma((i\oplus j)-i)=0$,

\item[(ii)] $\sigma(m+j)+\sigma((m\oplus j)+m)=1$.
\end{enumerate}
\end{lemma}
\begin{proof}
We have $(i\oplus j)-i = j-2m\sigma(i+j)$. Therefore $\sigma((i\oplus j)-i)$ is
opposite to $\sigma(i+j)$. This shows (i).

Let us prove (ii). If $j\le 0$, we have $\sigma(m+j)=0$ and $\sigma((m\oplus
j)+m)=\sigma(2m+j)$. Since $2m+j>2m-m=m$, we are done. If $j>0$, we have
$\sigma(m+j)=1$, and $\sigma((m\oplus j)+m)=\sigma(j)=0$.
\end{proof}

When $(G,*)$ is a obtained from $G$ by the cyclic or the dihedral construction,
denote by $x^*$ the inverse of $x$ in $(G,*)$, and by $x_i$ the $i$th power of
$x$ in $(G,*)$.

\begin{lemma}\label{Lm:Inverses} Let $G(*)=G[S,\alpha,h]$ be a cyclic
modification of $G$ such that $|G/S|=2m$. Then for $x\in G$ we have
\begin{displaymath}
     x^*=\left\{\begin{array}{ll}x^{-1},&x\not\in\alpha^m,\\
     x^{-1}h^{-1},&x\in\alpha^m.\end{array}\right.
\end{displaymath}
\end{lemma}
\begin{proof}
Assume that $x\in\alpha^i$, $i\in M\setminus\{m\}$. Then
$x^{-1}\in\alpha^{-i}$, and therefore
$x*x^{-1}=xx^{-1}h^{\sigma(0)}=xx^{-1}=1$. Assume that $x\in\alpha^m$. Then
$x^{-1}\in\alpha^{-m}=\alpha^m$. Therefore $x^{-1}h^{-1}\in\alpha^m$, too,
and we have $x*(x^{-1}h^{-1})=xx^{-1}h^{-1}h^{\sigma(m+m)} = 1$.
\end{proof}

\begin{lemma}\label{Lm:Conjugations} Under the assumptions of Lemma
$\ref{Lm:Inverses}$, we have $x*y*x^*=xyx^{-1}$, $x^**y*x=x^{-1}yx$ for every
$x$, $y\in G$.
\end{lemma}
\begin{proof} We only prove $x*y*x^*=xyx^{-1}$. The other equality is proved
along similar lines. Let $x\in\alpha^i$, $y\in\alpha^j$, $i\in M$, $j\in M$.

First assume that $i\ne m$. Then, by Lemma \ref{Lm:Inverses},
$x^*=x^{-1}\in\alpha^{-i}$, and we have
$x*y*x^*=xyh^{\sigma(i+j)}*x^{-1}=xyx^{-1}h^{\sigma(i+j)+\sigma((i\oplus
j)-i)}$. We are done by Lemma \ref{Lm:Aux}(i).

Now assume that $i=m$. Then, by Lemma \ref{Lm:Inverses},
$x^*=x^{-1}h^{-1}\in\alpha^m$, and we have $x*y*x^* =
xyh^{\sigma(m+j)}*(x^{-1}h^{-1})=xyx^{-1}h^{\sigma(m+j)+\sigma((m\oplus
j)+m)-1}$. We are done by Lemma \ref{Lm:Aux}(ii).
\end{proof}

\begin{lemma}\label{Lm:Powers} Assume that $(G,*)=G[S,\alpha,h]$, $|G/S|=2m$,
and $x\in\alpha$. Then
\begin{displaymath}
    x_i=\left\{\begin{array}{ll}x^i,&i\in M,\\ x^ih,&m<i\le
    2m.\end{array}\right.
\end{displaymath}
Furthermore, if $x\in\alpha^j$ and $j\in M$, we have $x_{2m}=x^{2m}h^j$.
\end{lemma}
\begin{proof}
First note that $x^i\in\alpha^i$ for every $i$. Therefore $x_i*x=x_ix$ for
every $i\in\{0,\dots,m-1\}$. This means that $x_i=x^i$ for every
$i\in\{0,\dots,m\}$.

Consider $x_i$ for $i\in\{-m+1,\dots,-1\}$. We have $(x_i)^*=x_{-i}$. By the
previous paragraph, $x_{-i}=x^{-i}=(x^i)^{-1}$. By Lemma \ref{Lm:Inverses},
$(x^i)^{-1}=(x^i)^*$. Altogether, we have $(x_i)^*=(x^i)^*$, and thus
$x_i=x^i$.

We have $x_m*x=x_mxh=x^mxh\in\alpha^{-m+1}$. It then follows that $x_i=x^ih$
for every $i\in\{m+1,\dots,2m\}$.

Let $x\in\alpha^j$, $j\in M$. Given $x$, $x'\in\alpha^k$, we have
$x_n=x^nh^\varepsilon$ and $x'_n=(x')^nh^\varepsilon$ for the same
$\varepsilon$, because the value of the exponent $\varepsilon$ depends only on
$n$ and $k$. We can therefore assume that $x=y^j$ for some $y\in\alpha$. Using
the above results, we have $x_{2m}=(y^j)_{2m}=(y_j)_{2m}=y_{2mj}=(y_{2m})_j =
(y_{2m})^j = (y^{2m}h)^j = (y^j)^{2m}h^j = x^{2m}h^j$.
\end{proof}

When $L$ is a Moufang loop, the associator subloop $A(L)$ can be defined
equivalently as the smallest normal subloop $H$ of $L$ such that $L/H$ is
associative. Therefore $A(L)\le S$ anytime $S$ is among the parameters of a
cyclic modification of $L$.

\begin{proposition}\label{Pr:Iso} Let $G_1=(G,\cdot)$, $G_2=(G,\circ)$
be two Moufang loops with common normal subloop $S$, and let $x\in G$ be such
that:
\begin{enumerate}
\item[(i)]  $G_1/S\cong G_2/S$ are cyclic of order $2m$,

\item[(ii)] both $G_1$, $G_2$ are generated by $S\cup\{x\}$,

\item[(iii)] the $2m$-th powers of $x$ coincide in $G_1$, $G_2$,

\item[(iv)] the conjugates $s^x$ for $s\in S$ coincide in $G_1$ and $G_2$,

\item[(v)] the multiplication in $S$ is the same in $G_1$, $G_2$,

\item[(vi)] the associators coincide in $G_1$, $G_2$,

\item[(vii)] $A(G_i)\le Z(G_i)\cap S$, for $i=1$, $2$,

\item[(viii)] $[a,b,cd]=[a,b,c\circ d]$ for every $a$, $b$, $c\in G$.
\end{enumerate}
Then $G_1$ is isomorphic to $G_2$.
\end{proposition}
\begin{proof}
Any element of $G_1$ decomposes uniquely as $x^is$, where $i\in
M=\{-m+1,\dots,m\}$, $s\in S$. Similarly, any element of $G_2$ decomposes
uniquely as $x_i\circ s$, where we use $x_i$ to denote the $i$th power of $x$
in $G_2$. Then the map $\varphi:G_2\to G_1$, $x_i\circ s\mapsto x^is$ is a
bijection.

We now show that $\varphi(x_k\circ s)=x^ks$ for every
$k\in\{-2m+2,\dots,2m\}$. When $k\in M$, we are done by the definition of
$\varphi$. Assume that $k>m$. Since $k-2m\in M$ and $x_{2m}$ is an element of
$S$, we have $\varphi(x_k\circ s)=\varphi(x_{k-2m}\circ x_{2m}\circ s) =
x^{k-2m}(x_{2m}\circ s)$. By (v), $x_{2m}\circ s = x_{2m}s$. Thus, by (iii),
$\varphi(x_k\circ s) = x^{k-2m}x^{2m}s=s$. Similarly, when $k<-m+1$, we have
$\varphi(x_k\circ s) = \varphi(x_{k+2m}\circ x_{-2m}\circ s) =
x^{k+2m}(x_{-2m}\circ s) = x^{k+2m}x_{-2m}s = x^{k+2m}x^{-2m}s=s$.

We also claim that $\varphi(s\circ x_k)=sx^k$ for $s\in S$,
$k\in\{-2m+2,\dots,2m\}$. Since $x_{-k}\circ s\circ x_k\in S$, we have
$\varphi(s\circ x_k) = \varphi(x_k\circ x_{-k}\circ s\circ x_k) =
x^k(x_{-k}\circ s\circ x_k)$. By (iv), the last expression is equal to
$x^ks$.

Define a new multiplication $*$ on $G$ by $x*y = \varphi(\varphi^{-1}(x)\circ
\varphi^{-1}(y))$. Then $(G,*)$ is isomorphic to $G_2$. We are going to show
that the multiplication $*$ coincides with the multiplication in $G_1$.

Now, for $i$, $j\in M$ and $s$, $t\in S$ we have $(x^is)*(tx^j) =
\varphi((x_i\circ s)\circ (t\circ x_j))$. By (vi) and (vii), $(x_i\circ
s)\circ (t\circ x_j) = x_i\circ(s\circ(t\circ x_j))\circ [x_i,s,t\circ x_j] =
x_i\circ(s\circ t)\circ x_j\circ [s,t,x_j]_{-1} \circ [x_i,s,t\circ x_j] =
x_i \circ (s\circ t)\circ [s,t,x_j]_{-1} \circ [x_i, s, t\circ x_j]\circ
x_{-i}\circ x_{i+j}$. By (v), (vii) and (viii), we can simplify this further
to $x_i((st)[s,t,x_j]^{-1}[x_i,s,tx_j])x_{-i} \circ x_{i+j}$. Therefore, by
the preceding paragraphs and (iv), $(x^is)*(tx^j) =
x^i(st)[s,t,x_j]^{-1}[x_i,s,tx_j]x^{-i}x^{i+j}$.

On the other hand, $(x^is)\cdot(tx^j) =
x^i(st)[s,t,x_j]^{-1}[x_i,s,tx_j]x^{-i}x^{i+j}$, and we are through.
\end{proof}

\begin{proposition}\label{Pr:Last} Suppose that $G$ is a Moufang $2$-loop such that
$A(G)\le Z(G)$. Suppose that $S$ is a normal subloop of $G$ such that $G/S$ is
cyclic of order $2m$, $G/S=\langle \alpha\rangle$. Let $j$ be relatively prime
to $2m$, and let $k\in M=\{-m+1,\dots,m\}$ be such that $jk\equiv 1\pmod{2m}$.
Then $G[S,\alpha^j,h]\cong G[S,\alpha,h^k]$.
\end{proposition}
\begin{proof}
Set $G_1=G[S,\alpha,h^k]$, $G_2=G[S,\alpha^j,h]$. Pick $x\in\alpha$. We are
going to check all assumptions of Proposition \ref{Pr:Iso}. By Lemma
\ref{Lm:Powers}, both $G_1$ and $G_2$ are generated by $S\cup\{x\}$, and the
$2m$-th power of $x$ in $G_1$ is equal to $x^{2m}h^k$. Since
$\alpha=(\alpha^j)^k$, the Lemma also implies that the $2m$-th power of $x$
in $G_2$ is equal to $x^{2m}h^k$. By Lemma \ref{Lm:Conjugations}, the
conjugates $s^x$ are the same in $G_1$ and $G_2$. The multiplication in $S$
is the same in $G_1$, $G_2$ (and $G$) by definition. By Theorem
\ref{Th:OldMain}, the associators of $G_i$ and $G$ are the same for $i=1$,
$2$. Thus the associators of $G_1$ and $G_2$ are the same. By the same
theorem, $A(G_i)\subseteq S$ for $i=1$, $2$. Consider the associators
$[a,b,cd]$, $[a,b,c\circ d]$, where $\cdot$ is the multiplication in $G$, and
$\circ$ is the multiplication in $G_1$. Note that $c\circ d$ differs from
$cd$ by a central element (namely a power of $h$). Therefore
$[a,b,cd]=[a,b,c\circ d]$. Similarly for $G$ and $G_2$. Thus all assumptions
of Proposition \ref{Pr:Iso} are satisfied, and $G_1\cong G_2$ follows.
\end{proof}

\begin{corollary} Let $G$ be a Moufang loop with normal subloop $S$ such that
$G/S$ is cyclic of order $2m$. Let $\alpha$, $\alpha'$ be two generators of
$G/S$. Then for every $h\in S\cap Z(G)$ there is $h'\in S\cap Z(G)$ such that
$G[S,\alpha,h]\cong G[S,\alpha',h']$.
\end{corollary}
\begin{proof} The generators of $G/S$ are exactly the powers $\alpha^j$, where
$(j,2m)=1$.
\end{proof}

\subsection{Speeding up step (iv)} The main bottleneck of the search is to
decide if the newly found loops $(L,*)$ are isomorphic to any of the
previously found loops. We are going to describe here how this problem was
overcome. In fact, it appears that the following algorithm performs very well
for all (power associative) loops, and $2$-loops in particular. Its idea is
natural and simple, but the details, based on theory and some heuristic, are
not so trivial.

Our task is to determine if two loops $L$, $M$ of order $n$ are isomorphic. The
main problem is that the space of possible isomorphisms is huge, consisting of
$n!$ bijections. Naturally, given an element $x$ of $L$, it cannot be mapped
onto an arbitrary element of $M$ if the mapping is supposed to be an
isomorphism. Certain invariants, such as the order of $x$, must be preserved.
The trick is to find invariants that are cheap yet powerful, in the sense that
the set of possible images of $x$ is small. Here are the invariants actually
used in the search:

For $x\in L$, let $I(x)=(|x|,s,f,(c_1,c_2,\dots,c_n))$, where
\begin{align*}
    s&=|\{y\in L;x=y^2\}|,\\
    f&=|\{y\in L;x=y^4\}|,\\
    c_i&=|\{y\in L; |y|=i,\,xy=yx\}.
\end{align*}
For a loop $L$ and an invariant $I$, let
\begin{align*}
    d_I&=|\{x\in L;\;I(x)=I\}|,\\
    D(L)& = \{ (I(x),d_{I(x)});\;x\in L\}.
\end{align*}
The distinguishing power of the \emph{discriminator} $D(L)$ is tremendous.
Table \ref{Tb:D} illustrates this eloquently for Moufang $2$-loops. For
instance, the table shows that the $2612$ loops forming the class $64:23$
give rise to $2331$ different discriminators in such a way that there are no
more than $6$ loops with the same discriminator. Hence, by precalculating the
discriminator once, at most $6$ instead of $2612$ loops have to be actually
tested for isomorphism at any given time in the search through the class
$64:23$.

Table \ref{Tb:D} also lists the maximum number of nonisomorphic modifications
$(L,*)$ of a loop $L$ in the given class. This shows that the constructions of
Section \ref{Sc:CD} often produce a large amount of nonisomorphic loops in one
step.

\begin{table}
\caption{The importance of discriminators in the search.}\label{Tb:D}
\begin{displaymath}
    \begin{array}{c|r|r|r|r}
    \text{class}&\text{$|$class$|$}&\text{discrim.\ types}&
    \text{max.\ with same discrim.}&\text{max.\ modifications}\\
    \hline
    16:01&5&5&1&3\\
    32:01&60&58&2&14\\
    32:04&11&11&1&5\\
    64:01&1316&1104&6&38\\
    64:04&18&18&1&6\\
    64:07&214&174&5&29\\
    64:15&11&11&1&5\\
    64:23&2612&2331&6&103\\
    64:38&44&44&1&19\\
    64:43&47&47&1&11
    \end{array}
\end{displaymath}
\end{table}


\section{The LOOPS package for GAP}\label{Sc:GAP}

\noindent The purpose of the GAP \cite{GAP} package LOOPS \cite{LOOPS} is to
implement calculation with loops and quasigroups in GAP. The package exists
only in a beta version and has not yet been accepted as a GAP shared package.
It is available online \cite{LOOPS}, together with installation instructions.

All Moufang loops found in this paper have now been included in the libraries
of LOOPS. Then $m$th nonassociative Moufang loop of order $n$ can be retrieved
by the command \texttt{MoufangLoop(n,m)}.

Since \cite{GoMaRa1999} already contains all nonassociative Moufang loops of
order less than $64$, LOOPS catalog numbers correspond to those of
\cite{GoMaRa1999}. Hence, for $n<64$, the Moufang loop called $n/m$ in
\cite{GoMaRa1999} is indeed isomorphic to \texttt{MoufangLoop(n,m)} of LOOPS.
The numbering of Moufang loops of order $64$ of LOOPS is based on our search.
For instance, the first 1316 Moufang loops of order $64$ are those of class
$64:01$.

Moreover, for a Moufang loop $L$ of order at most $64$, the LOOPS command
\texttt{IsomorphismTypeOfMoufangLoop(L)} returns the catalog number of $L$ and
the corresponding isomorphism, if possible. This command will be handy in the
search for additional Moufang loops of order $64$, should they exist.


\section{Conjectures}

\begin{conjecture} There are $4262$ nonassociative Moufang loops of order $64$,
as listed in this paper.
\end{conjecture}

The above conjecture holds if the following statement is true for $2^n=64$:
\emph{Every nonassociative Moufang $2$-loop of order $2^n$ is a modification
of a loop $M(G,2)$, where $G$ is a nonabelian group of order $2^{n-1}$.} In
view of Remark \ref{Rm:1}, the word ``nonassociative'' is essential in the
statement. Is the statement true for $64$? Is it true for $128$?

Finally, it is customary to classify loops with respect to isotopism in
addition to isomorphism. Recall that two loops $L$, $H$ are \emph{isotopic}
if there are bijections $\alpha$, $\beta$, $\gamma:L\to H$ such that
$\alpha(x)\beta(y)=\gamma(xy)$ for every $x$, $y\in L$. We ask: \emph{How do
the modifications behave with respect to isotopism? How many isotopism
classes of nonassociative Moufang loops of order $64$ are there?}


\section{Acknowledgement}

\noindent I would like to thank Edgar G.\ Goodaire for providing me with
electronic multiplication tables of Moufang loops of order at most $32$. I
also thank anonymous referees for several useful comments, and for an
improvement of the proof of Proposition \ref{Pr:Iso}.

\bibliographystyle{plain}

\end{document}